\documentclass{article}
\usepackage{amsmath, amssymb}

\usepackage[a4paper, margin=2cm]{geometry}

\title{GENERALIZED JENSEN'S INEQUALITY MOTIVATED FROM THERMODYNAMICS}
\author{Gourav Banerjee}
\date{20th October,2025}

\begin{document}


\maketitle

\begin{abstract}
In this paper, we generalize the work of P.T.Landsberg\cite{web1,web2} and S.S.Sidhu\cite{web3} by providing an inequality that has its main motivation from the laws of thermodynamics, in the form of a theorem which is quite useful in generating different inequalities such as the weighted AM-GM-HM inequality, the p-th power inequality , Jensen's inequality and many other inequalities.In this paper, we have not only given the thermodynamic motivation behind the inequality but we have given the required mathematical justification in the form of a straightforward rigorous proof using basic real analysis , which was not present in the works of Landsberg and Sidhu.
In fact, the first statement of the theorem mathematically proves the uniqueness of the equilibrium temperature that is attained when n different bodies at different temperatures are brought in contact.
The second statement of the theorem gives a mathematical proof of the fact that the  process in which n bodies at different temperatures  when brought in contact equilibriate to a common temperature is spontaneous,i.e., entropically favourable.

  Thus, this article motivates the students to come up with different mathematical results by observing the phenomena already existing in nature and also helps them to appreciate the conventional inequalities taught to them at the secondary school and undergraduate level by associating relevant physical phenomena with those inequalities.
\end{abstract}
\section*{Keywords} 
WEIGHTED AM-GM-HM INEQUALITY ; JENSEN INEQUALITY ; CONVEX FUNCTION ; THERMODYNAMICS

\vspace{.12 in}


\section*{Introduction}

One of the most commonly observed phenomenon in nature is that of systems coming to equilibrium.

One such phenomenon is that different bodies attain an equilibrium temperature when brought in contact at different temperatures.
This happens because it is an entropically favourable process(The net entropy of the system comprising of the bodies increases).
However, the different bodies may have widely varying thermal capacities which are a function of their temperatures, still the process of attaining an equilibrium temperature will be an entropically favourable process.

The theorem presented in this paper rigorously answers why the above mentioned process is entropically favourable and it also explains the existence of a unique equilibrium temperature attained by all the different bodies brought in contact with each other.

Students are introduced to many useful mathematical inequalities while pursuing undergraduate courses, such as the weighted Arithmetic mean-Geometric Mean- Harmonic mean inequality, pth-power inequality, etc.

The conventional proofs\cite{book:E,book:T} for the above-mentioned inequalities are lengthy and involve using Jensen's inequality for convex functions ( a real-valued function is called convex if the line segment between any two distinct points on the graph of the function lies above the graph between the two points. A twice-differentiable function of a single variable is convex if and only if its second derivative is non-negative on its entire domain) or some other computational trick. Also, the proofs do not generally encompass any connection of these inequalities to the natural world.

The theorem provided below which derives its main motivation from the laws of thermodynamics, helps us to easily generate all the above mentioned inequalities,along with Jensen's inequalities and even some other non-conventional inequalities without involving much computation.

\section*{The Main Result}
Let there be a set of $n$ arbitrary continuous functions,
\[ S = \{f_i\;|\; \mbox{where}\; f_i \in C((0,\infty))\;\; i \in \{1,2,\cdots,n\}\; \mbox{and}\; f_i : \mathbb{R}_{>0}  \rightarrow \mathbb{R}_{>0} \}\]

and let there be $n$ positive real numbers such that
$0<x_1<x_2<\cdots <x_n$
then

    1) There exists a unique $x_0$ such that
          $$
              \sum_{i=1}^{n}{\int_{x_i}^{x_0}{f_i}}  = 0
          $$
    2)  \(\displaystyle\sum_{i=1}^{n}{\int_{x_i}^{x_0}{f_ig}}  \geq 0\;\;\;\;\mbox{if g(x) is a decreasing function and integrable in } [x_1,x_n] \)         
              
    3)  \(\displaystyle\sum_{i=1}^{n}{\int_{x_i}^{x_0}{f_ig}}  \leq 0\;\;\;\;\mbox{if g(x) is an increasing function and integrable in } [x_1,x_n]\)
              
Before presenting the rigorous proof to this theorem, first let us present a short motivation from the laws of thermodynamics which has led me towards framing this theorem.

\section{MOTIVATION FROM THERMODYNAMICS}

 Note that this section just provides the inspiration from thermodynamics that led me to frame the above-mentioned inequality, and it is not be interpreted as the source of the generalization to Jensen's inequality.Rather, the explanation to the upcoming physical phenomenon can be treated as an application of the above mentioned inequality.
 
Consider a scenario where \textit{n} bodies indexed as $A_i$'s having thermal capacities equal to $f_i$'s(where $f_i$'s are continuous functions of temperature) and temperature equal to $x_i$'s are brought in contact and then insulated from the surroundings.

The Zeroth law \cite{web} of thermodynamics tells us that at equilibrium the temperature of all the bodies in contact will be the same and this temperature is unique (note that here we have assumed that none of the bodies undergo phase transition).

Now we apply the first law\cite{book:h} of thermodynamics which is energy conservation. Since the \textit{n} bodies had been isolated from surroundings, the net heat exchanged by the system should be 0. The net heat exchanged by the system is the sum of the amounts of heat exchanged by the individual bodies. Heat exchanged by the body $A_i$ is given by
\[
     \int_{x_i}^{x_0} f_i(x) \; dx
\]
therefore  the total heat exchanged by the system is
\[
    \sum_{i=1}^{n} \int_{x_i}^{x_0} f_i(x)\; dx = 0
\]
Thus we get the motivation behind the first statement of the theorem.

Since the phenomenon of heat exchange between the bodies is entropically favourable, therefore the total entropy change of the system during this process must be greater than equal to  0\cite{book:h}. The total entropy change of the system is the sum of the entropy change of the individual bodies present in system.The entropy change for body $A_i$ is given by
\[
    \int_{x_i}^{x_0} \frac{dQ_{rev}}{x} \;\;\;\;\;\; \mbox{where  } dQ_{rev} \mbox{  is the change in heat energy of a body which it undergoes reversibly }
\]
and $x$ is the instantaneous temperature of the body.
\[
    \int_{x_i}^{x_0} \frac{dQ_{rev}}{x} = \int_{x_i}^{x_0} \frac{f_i(x)}{x} \;dx
\]
Therefore the total entropy change of the system is given by:
\[
    \sum_{i=1}^{n} \int_{x_i}^{x_0} \frac{f_i(x)}{x} \;dx \geq 0
\]
Thus this is the motivation behind the second statement of the theorem and the function g(x)  in this case is $\frac{1}{x}$ which is a  decreasing function.
Thus this physical scenario of heat transfer between \textit{n} different bodies motivates us towards framing the statements of the above mentioned theorem. Now, I present the rigororous proof of the above mentioned theorem.

\section{PROOF:}

Consider the function F(x) = $\displaystyle
    \sum_{i=1}^{n}{\int_{x_i}^{x}{f_i}}  = 0
$\\
F(x) is lipschitz continous from second fundamental theorem of calculus.

Note that $F(x_1)\leq 0$ and $F(x_n)\geq 0$, hence by Intermediate value theorem , there exists
$\; x_0 \in [x_1,x_n]$
such that 

$$ F(x_0) = 0$$
let 
$x_h >x_l$
\begin{eqnarray*}
    \mbox{then}\;\; F(x_h)  = &\displaystyle  \sum_{i=1}^{n}{\int_{x_i}^{x_h}{f_i}}  = 0\\
    = & \displaystyle \sum_{i=1}^{n}{\left(\int_{x_i}^{x_l}{f_i} + \int_{x_l}^{x_h}{f_i}\right)}\\
    = & \displaystyle \sum_{i=1}^{n}{\int_{x_i}^{x_l}{f_i} + \sum_{i=1}^{n}{\int_{x_l}^{x_h}{f_i}}}\\
    & = F(x_l)+C \;\; \mbox{where}\;\; C>0
\end{eqnarray*}
hence $F(x_h)>F(x_l)$\\
so, $ F(x)$ is strictly increasing\\
Therefore the existance of $x_0$ is unique

Let $k \in \{0,1,\cdots ,n-1\}$ be such that $x_0\in [x_k,x_{k+1}]$ 

$\mbox{Since } g(x)\mbox{ is  decreasing}$
\begin{eqnarray}
    g(x_0) \leq g(x) \;\;\;\;\;\;\;\;\;\;\;\;\;\;\; &\;\;\;\;\;\;&  \mbox{for all} \;\;\; x   \leq  x_0  \\
    \int_{x_i}^{x_0}f_i(x)g(x_0)  \leq \int_{x_i}^{x_0}f_i(x)g(x)   &&\mbox{for all} \;\;\; x_i  \leq x_k
\end{eqnarray}
Similarly
\begin{eqnarray}
    g(x_0)  \geq  g(x) \;\;\;\;\;\;\;\;\;\;\;\;\;\;\; && \mbox{for all} \;\;\; x \geq x_0 \\
    \int_{x_i}^{x_0}f_i(x)g(x_0)  \leq  \int_{x_i}^{x_0}f_i(x)g(x) &&\mbox{for all} \;\;\; x_i \geq x_{k+1}
\end{eqnarray}
by adding all the $i$'s
we get
\begin{eqnarray*}
    \sum_{i=1}^{n}{\int_{x_i}^{x_0}{f_i(x)g(x_0)}}  & \leq & \sum_{i=1}^{n}{\int_{x_i}^{x_0}{f_i(x)g(x)}} \;\;\;\;\;\;\;\;\;\mbox{from (2),(4)}\\
    0 & \leq & \sum_{i=1}^{n}{\int_{x_i}^{x_0}{f_i(x)g(x)}}
\end{eqnarray*}
Similarly for $g(x)$ being increasing , we get

\begin{eqnarray*}
    \sum_{i=1}^{n}{\int_{x_i}^{x_0}{f_i(x)g(x_0)}}  & \geq & \sum_{i=1}^{n}{\int_{x_i}^{x_0}{f_i(x)g(x)}} \\
    0 & \geq & \sum_{i=1}^{n}{\int_{x_i}^{x_0}{f_i(x)g(x)}}
\end{eqnarray*}
\textbf{Note:}
While proving this inequality we have just used the property of integrability of $f_i$'s and monotonicity of $g$. in this proof we haven't used continuity of $f_i$'s thus this inequality also holds true if the functions $f_i$'s are non-negative integrable functions, however, the continuity property just gives us the uniqueness of $x_0$.

\section{APPLICATIONS}

\subsection{Proving Jensen's inequality for a differentiable function:}
Jensen's inequality states that:

$$F\left(\frac{\sum_{i=1}^{n}\lambda_ix_i}{\sum_{i=1}^{n} \lambda_i}\right) \leq \frac{\sum_{i=1}^{n}\lambda_iF(x_i)}{\sum_{i=1}^{n}\lambda_i}\;\;\; for\;\; \lambda_i's >0 \mbox{ if F is a convex function }$$

$$    F\left(\frac{\sum_{i=1}^{n}\lambda_ix_i}{\sum_{i=1}^{n} \lambda_i}\right) \geq \frac{\sum_{i=1}^{n}\lambda_iF(x_i)}{\sum_{i=1}^{n} \lambda_i}\;\;\; for\;\; \lambda_i's >0 \mbox{ if F is a concave function}$$

Proof:
For convex/concave functions ,note that F'(x) is monotonic,then take $f_i$'s to be $= \lambda_i$'s and $g(x) = F'(x)$ and the inequality follows straightway from the theorem.

$f_i = \lambda_i \mbox{ then } \displaystyle   \sum_{i=1}^{n}{\int_{x_i}^{x_0}{\lambda_i}} = \sum_{i=1}^{n} c_i(x_0 - x_i)  = 0 $\\
That implies,   $x_0 = \frac{(\lambda_1x_1 +\dots +\lambda_n x_n)}{(\lambda_1+\dots +\lambda_n)}$\\
$$    \sum_{i=1}^{n}\int_{x_i}^{x_0}{\lambda_iF'(x)} = \sum_{i=1}^{n}{\lambda_i(F(x_0) - F(x_i))}  $$
\\ Therefore 
$$F(x_0) \leq \frac{\sum_{i=1}^{n} \lambda_i F(x_i)}{\sum_{i=1}^{n}\lambda_i} \mbox{ if F is a convex function, } F(x_0) \geq \frac{\sum_{i=1}^{n} \lambda_i F(x_i)}{\sum_{i=1}^{n}\lambda_i} \mbox{ if F is a concave function. }$$

Thus inequalities which can be directly proved from Jensen' inequality can also be proved from the above-stated inequality.

We state some of those inequalities here which can be directly obtained from Jensen's inequality without explicitly proving them, but just by denoting the suitable $ f_i's \mbox{ and } g(x)$ that we need to consider :
\subsection{AM-GM-HM inequality}
 $\mbox{Consider} \;\;f_i = c_i\mbox{ and } \; g(x)=\frac{1}{x}\mbox{.}$ 
 Hence we have,
\[
    x_0 = \frac{\sum_{i=1}^{n} c_i x_i }{C_T} \;\;\;\;\;\mbox{ where } C_T \mbox{ is the sum of all } c_i's 
 \]\\
 $$    x_0 \geq \sqrt[C_T]{x_1^{c_1} x_2^{c_2} \cdots x_n^{c_n}} \;\;\;\;\mbox{ where } C_T \mbox{ is the sum of all } c_i's $$


 \subsection{$p^{th}$ power inequality}
 $    \mbox{consider } f_i = c_i\mbox{and } g(x)=x^{p-1} $ 
 , for $p \in (0,1)$ we find that $g(x)$ is  decreasing
 $$    \sum_{i=1}^{n}{\int_{x_i}^{x_0}{c_i x^{p-1} dx}}\geq 0 \implies  x_0^p \geq \frac{\sum_{i=1}^{n} c_ix_i^p}{C_T}$$
 likewise for $p \in (1,\infty)$
 $$    x_0^p \leq \frac{\sum_{i=1}^{n} c_ix_i^p} {C_T}  \;\;\;\mbox{where $C_T$ is the sum of all $c_i$'s}$$



\subsection{Some other inequalities}

Now we frame some other inequalities which cannot be directly proved from Jensen's inequality. Note that the key observation in proving Jensen's inequality was to choose the $f_i's$ to be constant functions, however if we choose our $f_i's$ to be some other functions then the corresponding inequalities that we will be genreating cannot be generated from Jensen's inequality.

$ GM(F(x)) $ : geometric mean of all the $F(x_i)'s$, $ AM(F(x))$ : arithmetic mean of all the $F(x_i)'s$ \\

$    \mbox{Now consider} \;\;\;f_i = \frac{c_i}{x}$
$\mbox{ let }  g(x) =  \frac{kx^k}{x^k+c}$
here $g$ is  increasing for all $x > 0$

$$   \sum_{i=1}^{n}c_i \int_{x_i}^{x_{GM}}  \frac{c_ikx^{k-1}}{x^k+c} \leq  0 \implies     (x_{GM}^{k} + c)^{C_T}  \leq  \prod_{i=1}^{n} (x_i^k + c)^{c_i}$$

Here $x_{GM}$ is the geometric mean of the $x_i$'s.



$    \mbox{Now let} \;\;\; g(x)= \frac{1-log(x)}{x}$
$    g'(x) < 0  $ for all $ x \in (0,e^2) , g'(x) > 0 $ for all $ x \in (e^2,\infty)$

$$    \sum_{i=1}^{n} c_i \int_{x_i}^{x_0} \frac{1-log(x)}{x^2} \geq 0  \;\;\;\;\; \mbox{ for all } x \in (0,e^2)$$
$$    C_T \left(\frac{log(x_0)}{x_0}\right)  \geq \sum_{i=1}^{n} c_i\frac{ log(x_i)}{x_i} \implies   \frac{log(GM(x))}{GM(x)}  \geq AM\left(\frac{log(x)}{x}\right)$$

as $\displaystyle\frac{1-log(x)}{x}$ is  decreasing in interval $(0,e^2)$ the inequality flips in interval $(e^2,\infty)$

$$\frac{log(GM(x))}{GM(x)}  \geq AM\left(\frac{log(x)}{x}\right) \;\;\;\;\;\; \mbox{ for all } x \in (0,e^2) $$
$$    \frac{log(GM(x))}{GM(x)}  \leq AM\left(\frac{log(x)}{x}\right) \;\;\;\;\;\; \mbox{ for all } x \in (e^2,\infty)$$
now for $f_i$ = $c_i$ and $g(x) = \frac{1-log(x)}{x^2}$
,$x_0$=AM ,$\;g'(x) < 0$ for all $x \in (0,e^{\frac{3}{2}})  ,  g'(x) > 0 $ for all $ x \in (e^{\frac{3}{2}},\infty)$

as $\displaystyle\frac{1-log(x)}{x^2}$ is  decreasing in interval $(0,e^\frac{3}{2})$,the inequality flips in interval $(e^\frac{3}{2},\infty)$

$$    \frac{log(AM)}{AM}  \geq AM\left(\frac{log(x)}{x}\right) \;\;\;\;\;\; \mbox{ for all } x \in (0,e^2) $$
$$    \frac{log(AM)}{AM}  \leq AM\left(\frac{log(x)}{x}\right) \;\;\;\;\;\; \mbox{ for all } x \in (e^2,\infty)$$

\section{CONCLUSION}

Thus we have successfully demonstrated that our theorem acts as a generator for Jensen inequalities as well as other inequalities. We just have to choose suitable $f_i 's $ and suitable $g(x)$ as our input functions for the theorem to easily give out these inequalities.

Thus, this paper might also serve as an impetus for students to derive other important mathematical inequalities like the Cauchy-Schwarz inequalities  from concepts in physics in the future.

\end{document}